\documentclass[a4paper,10pt]{article}
\usepackage{amsfonts,amssymb,euscript,amsmath,cite,graphicx}

\begin{document}

\title{\bf A New Version of a Posteriori Choosing \\
 Regularization Parameter in Ill-Posed Problems}

\author{\bf V. S. Sizikov}
\date{}
\maketitle

\begin{center}
SPbSU ITMO, Kronverksky pr 49, \\
St. Petersburg, 197101 Russia \\
e-mail: sizikov2000@mail.ru
\end{center}

\textbf{Abstract.} The new version of a posteriori choice (NVAC)
of the regularization parameter $\alpha$ in the classical
Tikhonov regularization method is considered. Lemmas and theorems
on the error and the asymptotic convergence rate of the regularized
solution are proved. A numerical example is given.



\bigskip

\textbf{1. Introduction}

\bigskip

Consider an operator equation of the first kind
\begin{equation}
\label{eq1}
A\,y = f, \quad y \in H_1, \quad f \in H_2 \,,
\end{equation}

\noindent
where $H_1$ and $H_2$ are Hilbert spaces and $A:\;H_1 \to H_2$
is a linear bounded operator. Suppose that the exact solution
$\bar{y}$ is the normal pseudosolution \cite{LeoYag,Siz}. Let, instead
of the exact $f$ and $A$, we have $\widetilde{f}$ and
$\widetilde{A}$ such that $\bigl\| \widetilde{f}-f\bigr\| \le \delta$,
$\delta > 0$, $\bigl\| \widetilde{A} - A\bigr\| \le \theta$, $\theta \ge 0$.
Denote by $\gamma \equiv (\delta, \theta)$. Given $\widetilde{f}$,
$\widetilde{A}$, $\delta$, and $\theta$, the problem is to find
an element $y_{\gamma} \in H_1$ that is a stable approximation
of $\bar{y}$ such that $\|y_{\gamma} - \bar{y}\| \to 0$ as
$\gamma \to 0$. \par

In the classical Tikhonov regularization method (using stabilizers of
the type $\|y\|_{L_2}^2$ or $\|y\|_{W_2^n}^2$), one solves the equation
\cite{LeoYag,Siz,VaiVer,VerSiz,BakGon1,BakGon2,TGSY,TLY,Koj,PetSiz}
\begin{equation}
\label{eq2}
\alpha\,y_{\alpha} + \widetilde{A}^* \widetilde{A}\,y_{\alpha}
= \widetilde{A}^* \widetilde{f} \,,
\end{equation}

\noindent
where $\alpha > 0$ is the regularization parameter. \par

Well-known \emph{ways for choosing the regularization parameter}
$\alpha$ were developed, namely, the discrepancy principle \cite{Mor},
the generalized discrepancy principle (GDP) \cite{TGSY}, the modified
discrepancy principle (MDP) \cite{Engl,Gfr,EngGfr,Gro1,GroScher,EHN},
the cross-validation method \cite{Golub}, the iteration stopping rule
by discrepancy \cite{BakGon1,BakGon2}, the local regularising algorithm
\cite{VosMukh}, the adaptive specialized generalized discrepancy
principle (SGDP) \cite{LeoYag}, etc. Estimates of
the error $\| y_{\alpha} - \bar{y}\|$ for the regularized
solution $y_{\alpha}$ were obtained, among them, with use
of an a priori information about the solution $\bar{y}$
(the sourcewise representability, etc.)
\cite{LeoYag,Siz,VaiVer,BakGon1,BakGon2,TGSY,TLY,Koj,Mor,Engl,Gfr,
EngGfr,Gro1,GroScher,EHN,Gro2,Zong}. \par

However, solving a number of model examples shows the following (see
\cite{VerSiz,TGSY}, et al.). For finite $\delta$ and $\theta$, the
principles can overstate the value of $\alpha$ in comparison with
$\alpha_{\rm opt}$. As a result, the error $\|y_{\alpha} - \bar {y}\|$
is overstated in comparison with $\|y_{\alpha_{\rm opt}} - \bar {y}\|$,
and the solution $y_{\alpha}$ becomes more smooth than
$y_{\alpha_{\rm opt}}$, and ``the fine structure'' of the solution
$y_{\alpha}$ is lost (cf. \cite{Vasin}).
Here, $\alpha_{\rm opt}$ is the value of $\alpha$ for which
$\|y_{\alpha} - \bar{y}\| = \min\limits_{\alpha}$
(the value of $\alpha_{\rm opt}$ can be determined without strong
a priori suppositions about the solution only in solving model
examples). This effect usually appears when the relative errors
$\delta_{\rm rel}$ and $\theta_{\rm rel} \gtrsim 1\%$
\cite[p.~283]{VerSiz}, \cite{TGSY}. \par

The aim of this paper is the further development of the new version
of a posteriori choice of $\alpha$ (NVAC) \cite{Siz} concentrating
attention on the question about closeness of $\alpha$ to
$\alpha_{\rm opt}$ and, as a result, of $\|y_{\alpha} - \bar{y}\|$
to $\|y_{\alpha_{\rm opt}} - \bar {y}\|$, furthermore, not so much
in asymptotics for $\delta, \theta \to 0$, as for finite $\delta$
and $\theta$. In this paper, the modified formulations of the
NVAC's statements are given, moreover, as far as possible without
using the sourcewise representability of $\bar {y}$.
In this case, the solution error estimates for finite $\delta$,
$\theta$, and $\alpha$ depend on the exact solution $\bar{y}$
that is known only in model examples. And in asymptotics
(for $\delta, \theta, \alpha \to 0$), the order of convergence
of $y_{\alpha}$ to $\bar{y}$ will be obtained.

\textbf{Remark 1.} Since $\alpha_{\rm opt}$ and $y_{\alpha_{\rm opt}}$
are known only in model examples but are unknown in real problems,
so the efficiency of the new version must be verified for model
examples.

\bigskip

\textbf{2. The idea of the NVAC}

\bigskip

Let us write Eq. (\ref{eq2}) in the form
\begin{equation}
\label{eq3}
\alpha\,y_{\alpha} + \widetilde{R}\,y_{\alpha} = \widetilde{F} \,,
\end{equation}

\noindent
where $\widetilde{R} = \widetilde{A}^* \widetilde{A}$,
$\widetilde{F} = \widetilde{A}^* \widetilde{f}$. \par

Along with the operator equation (\ref{eq1}), consider the Fredholm
integral equation of the first kind
\begin{equation}
\label{eq4}
A\,y \equiv \int_a^b K(x,s)\,y(s)\,ds = f(x),
\quad c \le x \le d \,.
\end{equation}

In the Tikhonov regularization method, instead of Eq. (\ref{eq4}),
one solves the equation (for $H_1=W_2^1$, $H_2=L_2$)
\cite[p.~24]{VerSiz}, \cite{TiAr}
\begin{equation}
\label{eq5}
\alpha \, [y_{\alpha}(t) - \tau \, y_{\alpha}''(t)] +
\int_a^b \widetilde{R}(t,s) \, y_{\alpha}(s) \, ds =
\widetilde{F}(t),
\quad a \le t \le b, \quad \tau \ge 0,
\end{equation}

$$
y_{\alpha}'(a) = y_{\alpha}'(b) = 0\,,
$$

\begin{equation}
\label{eq6}
\widetilde{R}(t,s) = \widetilde{R}(s,t) =
\int_c^d \widetilde{K}(x,t)\,\widetilde{K}(x,s)\,dx \,,
\end{equation}

\begin{equation}
\label{eq7}
\widetilde{F}(t) = \int_c^d \widetilde{K}(x,t)\,
\widetilde{f}(x)\,dx \,.
\end{equation}

Actually, the original equation in the Tikhonov regularization method
is the equation $\widetilde{A}^* \widetilde{A} \, y = \widetilde{A}^*
\,\widetilde{f}$ rather than $\widetilde{A}\,y = \widetilde{f}$.
Different variants of the discrepancy principle
\cite{VaiVer,VerSiz,BakGon1,BakGon2,TGSY,TLY,Koj,
Mor,Engl,Gfr,EngGfr,Gro1,GroScher,EHN,Gro2,Zong}
use the error $\delta$ of the right-hand side $\widetilde{f}$.
However, the function $\widetilde{f}(x)$ does not appear explicitly
as a right-hand side in the Tikhonov method. The right-hand side is
the function $\widetilde{F}(t)$ (see (\ref{eq3}) and (\ref{eq5})).
The function $\widetilde{f}(x)$ comes under the integral sign in the
expression for $\widetilde{F}(t)$ (see (\ref{eq7})), while the
integration operation is a smoothing filter with respect to
$\widetilde{f}(x)$. As a result, random errors in $\widetilde{f}(x)$
will be smoothed to a certain extent. In this case, the relative
error in $\widetilde{F}(t)$ can become considerably less than the
relative error in $\widetilde{f}(x)$ \cite{Siz}. \par

Concerning the error $\theta$ of the operator $\widetilde{A}$, the
factual operator in the Tikhonov method is the operator
$\widetilde{R} \equiv \widetilde{A}^* \widetilde{A}$ rather than
$\widetilde{A}$. Therefore, in choosing $\alpha$ from a discrepancy,
it is more appropriately to use the errors of the elements
$\widetilde{F}$ and $\widetilde{R}$ rather than $\delta$ and
$\theta$ (the errors of $\widetilde{f}$ and $\widetilde{A}$).
However, on deriving asymptotic estimates for $\alpha$ and for an
error of the solution $y_{\alpha}$, one should use the errors of
both the elements $\widetilde{F}$ and $\widetilde{R}$ and ones
$\widetilde{f}$ and $\widetilde{A}$.

In the generalized discrepancy principle (GDP) \cite{TGSY},
$\alpha=\alpha_{\rm d}$ (from discrepancy) is chosen to be a root
of the equation
$\bigl\| \widetilde{A}\,y_{\alpha} - \widetilde{f}\,\bigr\|^2 =
(\delta + \theta\,\|y_{\alpha}\|)^2 + \widetilde{\mu}^2$, where
$\widetilde{\mu} = \inf\limits_y\,\bigl\| \widetilde{A}\,y -
\widetilde{f}\,\bigr\|$ is the incompatibility measure of the equation
$\widetilde{A}\,y = \widetilde{f}$.

According to the Kojdecki way \cite{Koj},
$\alpha$ is a root of the equation
\begin{equation}
\label{eq8}
\alpha^q \, \bigl\| \widetilde{A}^* \widetilde{A} \, y_{\alpha}-
\widetilde{A}^* \widetilde{f} \, \bigr\| =  \beta \, \bigl\|
\widetilde{A} \bigr\| \left( \delta + \theta \, \| y_{\alpha} \| \right)
\end{equation}

\noindent
or, with regard to (\ref{eq2}),
$$
\alpha^{q+1} \| y_{\alpha} \| = \beta \, \bigl\| \widetilde{A} \bigr\|
\left( \delta + \theta \, \| y_{\alpha} \| \right) \,,
$$

\noindent
where $q \ge 0$ and $\beta > 0$ are some numbers.
One has proved \cite{Siz} the following lemma.

\textbf{Lemma 1.} \emph{The incompatibility measure}
$\widetilde{\nu} = \inf\limits_y \, \bigl\| \widetilde{R}\,y -
\widetilde{F}\bigr\|$ \emph{of the equation} $\widetilde{R}\,y =
\widetilde{F}$ \emph{is equal to zero.}

Now, we formulate again the new version of the a posteriori choice
of $\alpha$ (NVAC), moreover, the results obtained in \cite{Siz} will
be given without proofs. According to the NVAC, with regard to Lemma 1,
the regularization parameter $\alpha$ is chosen to be a root of the
equation \cite{Siz}
\begin{equation}
\label{eq9}
\alpha^q \, \bigl\| \widetilde{R}\,y_{\alpha} - \widetilde{F} \bigr\| =
\beta (\Delta + \Theta \, \|y_{\alpha}\|), \quad q \ge 0,
\quad \beta > 0 \,,
\end{equation}

\noindent
or a root of the equivalent equation
\begin{equation}
\label{eq10}
\alpha^{q+1} \|y_{\alpha}\| = \beta (\Delta + \Theta \, \|y_{\alpha}\|),
\quad q \ge 0, \quad \beta > 0 \,,
\end{equation}

\noindent
furthermore, $\bigl\| \widetilde{F} - F\bigr\| \le \Delta$ and $\bigl\|
\widetilde{R} - R \bigr\| \le \Theta$, where $\Delta = \Delta(\delta,\theta)
> 0$ is an upper estimate for the error of the right-hand side
$\widetilde{F}$ and $\Theta = \Theta(\theta) \ge 0$ is an upper estimate
for the error of the operator $\widetilde{R}$. Denote by
$\Gamma \equiv (\Delta,\Theta)$ and by $\alpha_{\rm n}$ a root of
(\ref{eq9}) or (\ref{eq10}) (the symbol ``n'' denotes ``new''). \par

\textbf{Remark 2.} Equation (\ref{eq9}) is rather like the equation
(\ref{eq8}). However, these equations have the difference of principle,
namely, in Eq. (\ref{eq8}), the errors $\delta$ and $\theta$ are used
and the factor $\bigl\| \widetilde{A}\,\|$ is separated from $\delta$ and
$\theta$, whereas in Eq. (\ref{eq9}), $\Delta$ and $\Theta$ are used.
Meanwhile, the value of $\bigl\| \widetilde{A}\,\bigr\| (\delta + \theta \,
\|y_{\alpha}\|)$ can be considerably greater than
$\Delta + \Theta \, \|y_{\alpha}\|$. This difference can lead to
overstated values of $\alpha$ and $\|y_{\alpha} - \bar{y}\|$. \par

\bigskip

\textbf{3. Justification of the New Version of a Posteriori Choosing
$\boldsymbol\alpha$}

\bigskip

Denote the left-hand side of (\ref{eq9}) or (\ref{eq10}) as
$$
\psi(\alpha) \equiv \alpha^q \bigl\| \widetilde{R} \, y_{\alpha} -
\widetilde{F} \, \bigr\| = \alpha^{q+1} \|y_{\alpha}\|
$$

\noindent
and the right-hand side of (\ref{eq9}) or (\ref{eq10}) as
$$
\xi(\alpha) \equiv \beta (\Delta + \Theta \, \|y_{\alpha}\|) \,.
$$

Then Eq. (\ref{eq9}) or (\ref{eq10}) can be written in the form
of the equation
\begin{equation}
\label{eq11}
\psi(\alpha) = \xi(\alpha) \,.
\end{equation}

\textbf{Lemma 2} \cite{Siz}. \emph{Under the condition}
\begin{equation}
\label{eq12}
\begin{array}{ll}
 \bigl\| \widetilde{F}\bigr\| > \beta\,\Delta, & q=0 \,, \\
 \bigl\| \widetilde{F}\bigr\| > 0, & q>0 \\
 \end{array}
\end{equation}

\noindent
\emph{the function} $\psi(\alpha)$ \emph{is continuous and
strictly monotonically increasing, moreover,}
$$
\lim\limits_{\alpha \to 0+} \, \psi(\alpha)=0 \,,
$$

$$
\lim\limits_{\alpha \to +\infty} \, \psi(\alpha) = \left\{
\begin{array}{cl}
 \bigl\| \widetilde{F}\bigr\|, & q=0 \,, \\
 0, & q>0 \quad and \quad \bigl\| \widetilde{F}\bigr\|=0 \,, \\
 \infty, & q>0 \quad and \quad \bigl\| \widetilde{F}\bigr\|>0 \,,
 \end{array} \right.
$$

\noindent
\emph{and function} $\xi(\alpha)$ \emph{is continuous and
strictly monotonically decreasing, moreover,}
$$
\begin{array}{ccl}
\lim\limits_{\alpha \to 0+}\,\xi(\alpha) & > &
\beta\,\Delta > 0 \,, \\
\lim\limits_{\alpha \to +\infty}\,\xi(\alpha) & = &
\beta\,\Delta > 0 \,.
\end{array}
$$

Now, the NVAC can be formulated as the following theorem.

\textbf{Theorem 1.} \emph{Let the equation} $\widetilde{A}
\, y = \widetilde{f}$, $y \in H_1$, $\widetilde{f} \in H_2$,
\emph{be solved by the Tikhonov regularization method according
to} (\ref{eq2}) \emph{or} (\ref{eq3}), \emph{where}
$\bigl\| \widetilde{f} - f\bigr\| \le \delta$, $\delta > 0$,
$\bigl\| \widetilde{A} - A\bigr\| \le \theta$, $\theta \ge 0$.
\emph{Suppose that the regularization parameter} $\alpha$
\emph{is chosen to be a root of Eq.} (\ref{eq9}), (\ref{eq10})
\emph{or} (\ref{eq11}), \emph{furthermore},
$\bigl\| \widetilde{F} - F \bigr\| \le \Delta$,
$\bigl\| \widetilde{R} - R \bigr\| \le \Theta$, \emph{where}
$\Delta = \Delta(\delta,\theta) > 0$, $\Theta = \Theta(\theta) \ge 0$.
\emph{Then, under condition} (\ref{eq12}), \emph{a root}
$\alpha = \alpha_{\rm n}$ \emph{of Eq.} (\ref{eq11})
\emph{exists and is unique, and the solution}
$y_{\alpha_{\rm n}}$ \emph{can be found by solving Eq.}
(\ref{eq3}) \emph{with} $\alpha = \alpha_{\rm n}$.
\emph{If condition} (\ref{eq12}) \emph{is not fulfilled,
then} $y_{\alpha_{\rm n}} = 0$.

\bigskip

\textbf{4.} \textbf{Some dependences}

\bigskip

Let us establish the dependences
$\Delta = \Delta(\delta,\theta)$ and $\Theta = \Theta(\theta)$.
The estimate for the error $\Delta$
of the right-hand side $\widetilde{F}$ has the form \cite{Siz}
\begin{equation}
\label{eq13}
\Delta \le \bigl\| \widetilde{A} \bigr\| \, \delta + \bigl\| \widetilde{f}
\bigr\| \, \theta \,,
\end{equation}
and the estimate for the error $\Theta$
of the operator $\widetilde{R}$ has the form \cite{Siz}
\begin{equation}
\label{eq14}
\Theta \le 2 \, \bigl\| \widetilde{A}\bigr\| \, \theta \,.
\end{equation}

\textbf{Remark 3.} The estimates (\ref{eq13}) and (\ref{eq14}) are
necessary for justifying the convergence of the NVAC. However, in
practice for a finite $\delta$ and $\theta$, the formulas (\ref{eq13})
and (\ref{eq14}) may give an overstatement of $\Delta$ and $\Theta$
(see example in the end of the present paper)
and, hence, of $\alpha_{\rm n}$ if one uses the upper estimates:
$\Delta = \bigl\| \widetilde{A}\bigr\| \, \delta + \bigl\| \widetilde{f}
\bigr\| \, \theta$ and $\Theta = 2 \, \bigl\| \widetilde{A} \bigr\| \, \theta$.
This overstatement is caused by that the factor $\bigl\| \widetilde{A}
\bigr\|$ is separated from $\delta$ and $\theta$ in the estimates
(\ref{eq13}) and (\ref{eq14}). To obtain more exact estimates of
$\Delta$ and $\Theta$, one can use, for example, the algorithms II,
III and V from the paper \cite{Siz}.

\bigskip

\textbf{5. Estimates for} $\boldsymbol\alpha_{\bf n}$

\bigskip

We give two upper estimates for $\alpha_{\rm n}$ in the NVAC.
Define \cite{Siz}, \cite[p.~78]{Koj}
\begin{equation}
\label{eq15}
\alpha_0 = \bigl\| \widetilde{R} \bigr\| = \bigl\| \widetilde{A} \bigr\|^2
= \bigl\| \widetilde{A}^* \bigr\|^2 \,.
\end{equation}

The condition (\ref{eq12}) for $q = 0$ can be written as
\begin{equation}
\label{eq16}
\frac{\Delta}{\bigl\| \widetilde{F} \bigr\|} < \frac1{\beta} \;.
\end{equation}

Let us introduce as an extended variant of condition (\ref{eq16})
the following condition \cite{Siz}
\begin{equation}
\label{eq17}
\frac{\Delta}{\bigl\| \widetilde{F} \bigr\|} + \frac{\Theta}{\bigl\|
\widetilde{R} \bigr\|} \le \frac1{\beta} \, \frac{\bigl\| \widetilde{R}
\bigr\|^q}2 \,.
\end{equation}

Condition (\ref{eq17}) can also be considered as a modification
of condition (53) in \cite{Koj}. It is proved \cite{Siz}

\textbf{Lemma 3.} \emph{Under condition} (\ref{eq17}),
\emph{one has the inequality}
\begin{equation}
\label{eq18}
\psi(\alpha_0) \ge \xi(\alpha_0) \,.
\end{equation}

\textbf{Corollary 1} \cite{Siz}. Since the functions $\psi(\alpha)$
and $\xi(\alpha)$ are increasing and decreasing, respectively,
relations (\ref{eq15}), (\ref{eq17}), (\ref{eq18}) imply that
\begin{equation}
\label{eq19}
\alpha_{\rm n} \le \alpha_0 = \bigl\| \widetilde{R}\bigr\| \,.
\end{equation}

Inequality (\ref{eq19}) gives an upper estimate for $\alpha_{\rm n}$
in terms of the norm of the operator. It is also proved \cite{Siz}

\textbf{Lemma 4.} \emph{Under condition} (\ref{eq12}),
\emph{it holds that}
\begin{equation}
\label{eq20}
\alpha_{\rm n} \le \left[ \beta \left( \frac{2\,\bigl\| \widetilde{R}
\bigr\|}{\bigl\| \widetilde{F} \bigr\|} \, \Delta + \Theta \right)
\right]^{1/(q+1)} \,.
\end{equation}

Inequality (\ref{eq20}) gives another upper estimate for
$\alpha_{\rm n}$ (in terms of the errors in the original data). \par

\textbf{Corollary 2} \cite{Siz}. Since
$$
\frac{2\,\bigl\| \widetilde{R}\bigr\|}{\bigl\| \widetilde{F}\bigr\|} \,
\Delta + \Theta \le \max \left\{ \frac{2\,\bigl\| \widetilde{R} \bigr\|}
{\bigl\| \widetilde{F} \bigr\|}, \, 1 \right\} \, (\Delta + \Theta) \,,
$$

\noindent
the estimate (\ref{eq20}) can be written as
\begin{equation}
\label{eq21}
\alpha_{\rm n} \le c_1 (\Delta + \Theta)^{1/(q + 1)} \,,
\end{equation}

\begin{equation}
\label{eq22}
c_1 = \left[ \beta \cdot \max \left\{ 2 \, \| \widetilde{R} \| /
\| \widetilde{F} \|, \, 1 \right\} \right]^{1/(q+1)}
> 0 \,.
\end{equation}

\textbf{Corollary 3} \cite{Siz}. Inequality (\ref{eq21}) generates the
asymptotic estimate
\begin{equation}
\label{eq23}
\alpha_{\rm n} = O \left( (\Delta + \Theta)^{1/(q+1)} \right) \,,
\quad \Delta, \Theta \to 0 \,.
\end{equation}

Using (\ref{eq13}) and (\ref{eq14}), we can write the estimates
(\ref{eq21}) and (\ref{eq22}) also as
\begin{equation}
\label{eq24}
\alpha_{\rm n} \le c_2 (\delta + \theta)^{1/(q + 1)} \,,
\end{equation}

\begin{equation}
\label{eq25}
c_2 = \left[2 \beta \, \| \widetilde{A} \| \cdot \max \left\{
\| \widetilde{R} \| / \| \widetilde{F} \|, \, \| \widetilde{A} \| \cdot
\| \widetilde{f} \| / \| \widetilde{F} \| + 1 \right\} \right]^{1/(q+1)}
> 0 \,,
\end{equation}

\begin{equation}
\label{eq26}
\alpha_{\rm n} = O \left( (\delta + \theta)^{1/(q+1)} \right) \,,
\quad \delta, \theta \to 0 \,.
\end{equation}

The relations (\ref{eq21}), (\ref{eq22}), (\ref{eq24}), (\ref{eq25})
show that the estimate for $\alpha_{\rm n}$ decreases with decrease
of $\beta$.

\bigskip

\textbf{6. Error Estimate for the Regularized Solution}

\bigskip

We give a new, more precise, estimate for the error
$\|y_{\alpha_{\rm n}} - \bar{y}\|$
of the regularized solution $y_{\alpha_{\rm n}}$ in the NVAC.
In the papers \cite{VaiVer,BakGon1,BakGon2,Koj,GroScher,Zong} et al.,
it was shown that in the Tikhonov regularization method
there holds the following error estimate for the regularized solution
(on the assumption that the exact solution $\bar{y}$ is sourcewise
representable with index 1, i.e. $\bar{y} = A^* A\,w$, $w \in H_1$):
\begin{equation}
\label{eq27}
\|y_{\alpha} - \bar{y}\| \le c_3 \, \frac{\delta + \theta}
{\sqrt{\alpha}} + c_4 \alpha \,,
\end{equation}
where $c_3, \, c_4 > 0$ are some constants.

Let us use the estimate (\ref{eq27}). For
$\alpha_{\rm n} = O \left( (\delta + \theta)^{1/(q+1)} \right)$
(see (\ref{eq26})) there exist such positive constants $a_1$ and $a_2$
that (cf. \cite[p. 65]{Koj})
\begin{equation}
\label{eq28}
a_1 (\delta + \theta)^{1/(q+1)} < \alpha_{\rm n} <
a_2 (\delta + \theta)^{1/(q+1)} \,.
\end{equation}

Hence,
\begin{equation}
\label{eq29}
\| y_{\alpha_{\rm n}} - \bar{y} \| \leq \frac{c_3}{\sqrt{a_1}} \,
(\delta + \theta)^{(q+0.5)/(q+1)} + c_4 a_2 \,
(\delta + \theta)^{1/(q+1)} \,.
\end{equation}

The estimate (\ref{eq29}) makes possible to obtain the following
\emph{asymptotic estimates}.

For sufficiently small $\delta$ and $\theta$, we have:
\begin{equation}
\label{eq30}
\| y_{\alpha_{\rm n}} - \bar{y} \| \le
c \, (\delta + \theta)^{\widetilde{q}} \,,
\quad c>0 \,,
\end{equation}

\begin{equation}
\label{eq31}
\widetilde{q} = \frac {\min\{q+0.5,\,1\}} {q+1} =
\left\{
\begin{array}{ll}
(q+0.5)/(q+1), \; & q \in [0,\,0.5]\,, \\
1/(q+1), \; & q \ge 0.5\,.
\end{array}
\right.
\end{equation}

As $\delta, \theta \to 0$, we obtain the asymptotic estimate
for the convergence rate of $y_{\alpha_{\rm n}}$ to $\bar{y}$:
\begin{equation}
\label{eq32}
\| y_{\alpha_{\rm n}} - \bar{y} \| = O \left( (\delta +
\theta)^{\widetilde{q}} \right),
\end{equation}

\noindent
as well as (we write again the estimate for $y_{\alpha_{\rm n}})$
\begin{equation}
\label{eq33}
\alpha_{\rm n} = O \left( (\delta + \theta)^{1/(q+1)} \right).
\end{equation}

The best asymptotic estimates are obtained for $q=0.5$:
\begin{equation}
\label{eq34}
\| y_{\alpha_{\rm n}} - \bar{y} \| = O \left( (\delta +
\theta)^{2/3} \right), \quad
\alpha_{\rm n} = O \left( (\delta + \theta)^{2/3} \right),
\end{equation}
i.e. the optimal order of convergence is obtained.
This is conform to results of the papers
\cite{Engl,Gfr,EngGfr,Gro1,GroScher,Zong} et al., in which
the optimal order of convergence has also been obtained, but
for other ways for choosing $\alpha$ (the modified discrepancy
principle, etc.).

If, e.g., $q=0$ then
$\| y_{\alpha_{\rm n}} - \bar{y} \| =
O \left( (\delta + \theta)^{1/2} \right)$ -- the suboptimal
order of convergence as in the GDP \cite{TGSY}.

\bigskip

\textbf{7. Final Theorem}

\bigskip

In conclusion, we prove the summarizing theorem.

\textbf{Theorem 2.} \emph{Let the equation} (\ref{eq2}) \emph{be
solved. Furthermore, the regularization parameter} $\alpha$
\emph{is chosen with the help of the NVAC according to} (\ref{eq11})
\emph{by equal} $\alpha = \alpha_{\rm n}$.
\emph{In this case, the estimates} (\ref{eq19})--(\ref{eq26})
\emph{for} $\alpha_{\rm n}$ \emph{and the estimates}
(\ref{eq29})--(\ref{eq32}) \emph{for the error}
$\|y_{\alpha_{\rm n}} - \bar{y}\|$
\emph{of the regularized solution} $y_{\alpha_{\rm n}}$
\emph{are valid. One has a convergence of the regularized solution}
$y_{\alpha_{\rm n}}$ \emph{to the exact solution} $\bar{y}$ \emph{as}
$\delta, \theta \to 0$, \emph{i.e. the NVAC generates a
regularizing algorithm.}

\textbf{Proof.} According to (\ref{eq30}), (\ref{eq32}),
$\| y_{\alpha_{\rm n}} - \bar{y} \| \to 0$ as $\delta, \theta \to 0$.
This means that
$y_{\alpha_{\rm n}}\xrightarrow[\delta, \theta \to 0]{} \bar{y}$.
Theorem 2 is proved.

\bigskip

\textbf{8. Numerical example}

\bigskip

To realize the new version of the a posteriori choice of $\alpha$, we
have developed the program package NVAC using Fortran PowerStation 4.0.
The following \emph{model example} (cf. \cite[p.~162]{PetSiz}) was
solved with the help of this package. \par

The exact solution was set as a superposition of five gaussians
(the solution with variations):
$$
\begin{array}{l}
\bar{y}(s) = 6.5\,e^{-[(s + 0.66)/0.085]^2} +  9\,e^{-[(s + 0.41)/0.075]^2} \\
            + 12\,e^{-[(s - 0.14)/0.084]^2} + 14\,e^{-[(s - 0.41)/0.095]^2}
             + 9\,e^{-[(s - 0.67)/0.065]^2} ,
\end{array}
$$

\begin{figure}[h]
\centering \includegraphics[height=6cm]{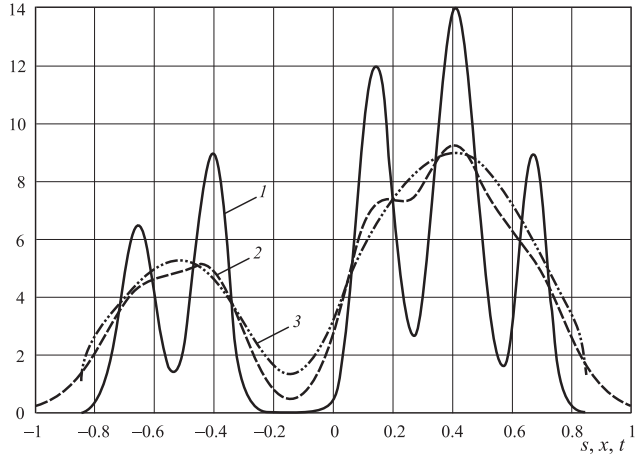}
\caption{\emph{1} --- $\bar{y}(s)$;
         \emph{2} --- $f(x)$;
         \emph{3} --- $F(t)$}
\label{f1}
\end{figure}

\noindent
$a = - 0.85$, $b = 0.85$, $c = - 1$, $d = 1$, the kernel
$$
K(x,s) = \sqrt{r/\pi}\,e^{-r(x - s)^2/(1 + x^2)} \,,
$$

\noindent
where the exact value $r$ is $r = 59.924$. The numbers of
discretization nodes are $l = 161$ (on $x$) and $n = 137$ (on $s$ and $t$).
The discretization steps are
$\Delta x = \Delta s = \Delta t = {\rm const} = 0.0125$.
In this example, $\|\bar{y}\| = 7.606$, $\|f\| = 6.907$, $\|A\| = 2.419$,
$\|F\|= 7.216$, $\|R\| = 2.196$. Figure \ref{f1} shows the exact solution
$\bar{y}(s)$, the right-hand side $f(x)$ (considerably more smooth than
$\bar{y}(s)$), and the new right-hand side $F(t)$ (still more smooth than
$f(x)$).

At first, the \emph{direct problem} was solved. The values $f_{i}$,
$i = 1,\ldots,l$, were calculated. The errors $\delta f_i$ distributed by
the normal law with zero expectation and with the mean square deviation
$\delta = 0.0001$, 0.15 and 0.5 were added to the values $f_i$.
The values $\widetilde{r} = 59.920$, 60 and 65 were used instead of
the exact value of $r$. Table 1 shows, as an instance, the values
of $\delta$, $\delta/\|f\|$, $\Delta=\|\Delta F\|$, $\Delta/\|F\|$ and
(for comparison) $\|\widetilde{A}\|\delta+\|\widetilde{f}\|\theta$ for
$\widetilde{r} = 60$. Such value of $\widetilde{r}$ corresponds to the
following parameters: $\theta=\|\Delta A\|=1.321 \cdot 10^{-3}$,
$\theta/\|A\| = 5.46 \cdot 10^{-4} = 0.0546\,\%$,
$\Theta=\|\Delta R\|=1.194 \cdot 10^{-3}$,
$\Theta/\|R\| = 5.44 \cdot 10^{-4} = 0.0544\,\%$,
$2\|\widetilde{A}\|\theta=6.392 \cdot 10^{-3}$.

\begin{center}
Table 1

\begin{tabular}{|c|c|c|c|c|}
\hline
$\delta$ & $\delta/\|f\|$ & $\Delta=\|\Delta F\|$ &
$\Delta/\|F\|$ & $\|\widetilde{A}\|\delta+\|\widetilde{f}\|\theta$ \\
\hline
0.0001 &
$\aligned 1.448 \cdot 10^{-5} \\ \approx 1.4 \cdot 10^{-3}\% \endaligned$ &
$0.6691 \cdot 10^{-3}$ &
$\aligned 0.927 \cdot 10^{-4} \\ \approx 0.93 \cdot 10^{-2}\% \endaligned$ &
$9.4 \cdot 10^{-3}$ \\
\hline
0.15 &
$\aligned 2.172 \cdot 10^{-2} \\ \approx 2.2\% \endaligned$ &
0.01878 &
$\aligned 0.259 \cdot 10^{-2} \\ \approx 0.26 \% \endaligned$ &
0.3721 \\
\hline
0.5 &
$\aligned 7.239 \cdot 10^{-2} \\ \approx 7.2\% \endaligned$ &
0.06256 &
$\aligned 0.867 \cdot 10^{-2} \\ \approx 0.87\% \endaligned$ &
1.219 \\
\hline
\end{tabular}
\end{center}

\smallskip

Furthermore, the operator norms  $\|A\|$,
$\theta = \bigl\| \widetilde{A} - A \bigr\|$, $\|R\|$, and
$\Theta = \bigl\| \widetilde{R} - R \bigr\|$ were calculated
by means of the Hilbert--Schmidt norm, e.g.,
$$
\|A\| = \left\{ \int_a^b \int_c^d K^2(x,s)\,dx\,ds \right\}^{1/2}.
$$

Comparing the values of $\Delta$ and
$\|\widetilde{A}\|\delta+\|\widetilde{f}\|\theta$, as well as
$\Theta$ and $2\|\widetilde{A}\|\theta$ (see (\ref{eq13}) and
(\ref{eq14})) we see that the upper estimates
$\|\widetilde{A}\|\delta+\|\widetilde{f}\|\theta$ and
$2\|\widetilde{A}\|\theta$ overstate by one order the values of
$\Delta$ and $\Theta$, and comparison of $\delta/\|f\|$ and
$\Delta/\|F\|$ shows that $\Delta/\|F\|$ less by one order than
$\delta/\|f\|$ for $\delta/\|f\| \gtrsim 1\%$.
About this, one says already above.

\begin{figure}[h]
\centering \includegraphics[height=6cm]{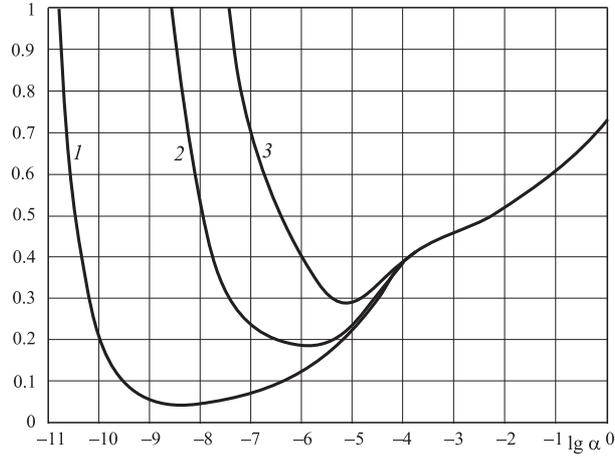}
\caption{The relative solution error
$\|y_{\alpha} - \bar{y}\| / \|\bar{y}\|$ at $\tau=1$
\emph{1} --- $\delta=0.0001$, $\widetilde{r}=59.920$;
\emph{2} --- $\delta=0.15$, $\widetilde{r}=60$;
\emph{3} --- $\delta=0.5$, $\widetilde{r}=65$}
\label{f2}
\end{figure}

Aftewards, the \emph{inverse problem} was solved. Equation (\ref{eq5})
was solved by the quadrature method at $\tau=1$ \cite[pp. 249--251]{VerSiz}.
Figure \ref{f2} shows some curves of the relative solution error
$\| y_{\alpha} - \bar{y} \| / \| \bar{y} \|$ (it can be calculated
only in solving a model example with known $\bar{y}$).

Table 2 shows, as an instance, the values of $\alpha_{\rm opt}$,
$\alpha_{\rm n}$ and the relative errors of the solutions
$y_{\alpha_{\rm opt}}$ and $y_{\alpha_{\rm n}}$ for $\widetilde{r}=60$,
$q=0$, $\tau=1$, $\beta=1$ and $\beta=0.1$\,.

\begin{center}
Table 2

\begin{tabular}{|c|c|c|c|c|c|c|}
\hline
$\delta$ & $\lg\alpha_{\rm opt}$ &
$\|y_{\alpha_{\rm opt}}-\bar{y}\| / \|\bar{y}\|$ &
\multicolumn{2}{|c|}{$\lg\alpha_{\rm n}$} &
\multicolumn{2}{c|}{$\|y_{\alpha_{\rm n}}-\bar{y}\| / \|\bar{y}\|$} \\
\cline{4-7}
& & & $\beta=1$ & $\beta=0.1$ & $\beta=1$ & $\beta=0.1$ \\
\hline
0.0001 & $-8.7$ & 0.0385 & $-5.1$ & $-6.2$ & 0.2107 & 0.1099 \\
\hline
0.15 & $-5.8$ & 0.1848 & $-4.3$ & $-5.7$ & 0.3466 & 0.1858 \\
\hline
0.5 & $-5.2$ & 0.2644 & $-3.6$ & $-5.2$ & 0.4311 & 0.2644 \\
\hline
\end{tabular}
\end{center}

\begin{figure}[p]
\centering \includegraphics[height=6cm]{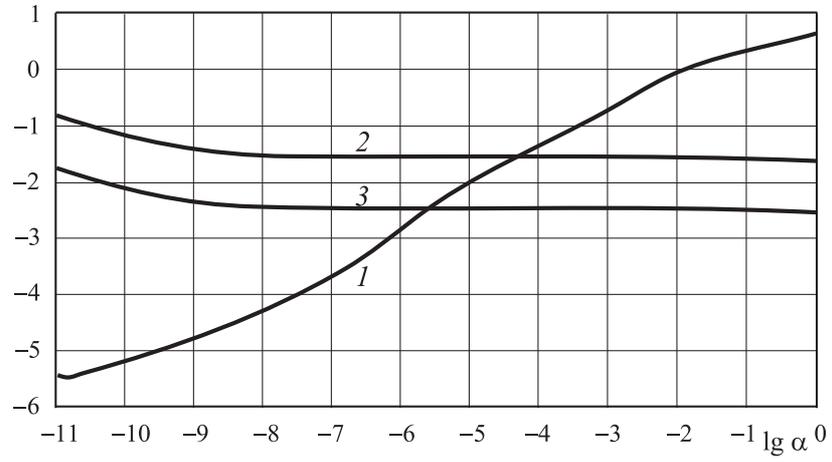}
\caption{\emph{1} --- $\lg\psi(\alpha)$;
         \emph{2} --- $\lg\xi(\alpha)$, $\beta=1$;
         \emph{3} --- $\lg\xi(\alpha)$, $\beta=0.1$}
\label{f3}
\end{figure}

Figure \ref{f3} shows the logarithms of the functions
$\psi(\alpha) = \alpha^q \bigl\| \widetilde{R}\,y_{\alpha}
- \widetilde{F} \bigr\|$ and $\xi(\alpha) = \beta \, (\Delta + \Theta \,
\| y_{\alpha} \|)$.

\begin{figure}[p]
\centering \includegraphics[height=8cm]{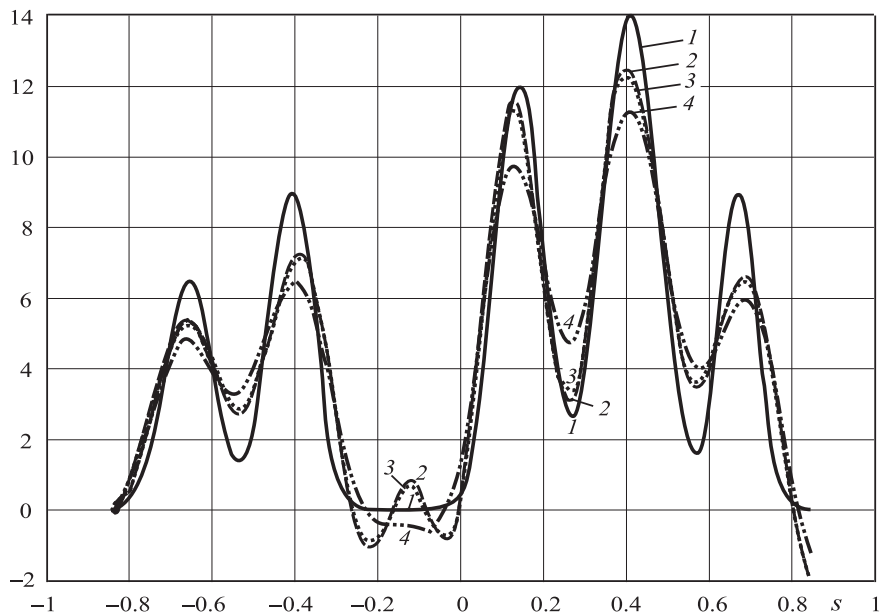}
\caption{\emph{1} --- $\bar{y}(s)$;
\emph{2} --- $y_{\alpha}(s)$, $\alpha=\alpha_{\rm opt}=10^{-5.8}$;
\emph{3} --- $y_{\alpha}(s)$, $\alpha=\alpha_{\rm n}=10^{-5.7}$; $\beta=0.1$;
\emph{4} --- $y_{\alpha}(s)$, $\alpha=\alpha_{\rm n}=10^{-4.3}$; $\beta=1$}
\label{f4}
\end{figure}

Figure \ref{f4} shows the exact solution $\bar{y}(s)$ and
the regularized solutions $y_{\alpha}(s)$ at
$\alpha = \alpha_{\rm opt} = 10^{-5.8}$,
$\alpha = \alpha_{\rm n} = 10^{-5.7}$ ($\beta=0.1$) and
$\alpha = \alpha_{\rm n} = 10^{-4.3}$ ($\beta=1$) for
$\delta=0.15$, $\widetilde{r}=60$, $q=0$, $\tau=1$.

\smallskip
The author thanks Prof. M.A. Kojdecki for useful discussion of paper
results.

\end{document}